\documentclass[11pt]{article}
\usepackage{amssymb,latexsym}
\usepackage{epsfig}
\usepackage{eufrak}
\usepackage{amsmath}
\usepackage{mathrsfs}
\usepackage{color}

\newtheorem{theorem}{Theorem}
\newtheorem{lemma}{Lemma}

\newcommand{\be}{\begin{equation}}
\newcommand{\ee}{\end{equation}}
\newcommand{\bee}{\begin{eqnarray*}}
\newcommand{\eee}{\end{eqnarray*}}
\newcommand{\bel}{\begin{eqnarray}}
\newcommand{\eel}{\end{eqnarray}}
\newcommand{\bec}{\begin{cases}}
\newcommand{\eec}{\end{cases}}
\newcommand{\bem}{\begin{bmatrix}}
\newcommand{\eem}{\end{bmatrix}}

\newcommand{\bed}{\begin{description}}
\newcommand{\eed}{\end{description}}
\newcommand{\bei}{\begin{itemize}}
\newcommand{\eei}{\end{itemize}}
\newcommand{\ben}{\begin{enumerate}}
\newcommand{\een}{\end{enumerate}}

\newcommand{\beL}{\begin{lemma}}
\newcommand{\eeL}{\end{lemma}}
\newcommand{\beT}{\begin{theorem}}
\newcommand{\eeT}{\end{theorem}}

\newcommand{\bpf}{\begin{pf}}
\newcommand{\epf}{\end{pf}}

\newcommand{\pfbox}{\hfill\mbox{$\Box$}}

\newenvironment{pf}{\paragraph*{Proof{\rm.}}}{\pfbox\bigskip}

\begin{document}

\title{{\bf On the Binomial Confidence Interval and Probabilistic Robust Control}}

\author{Xinjia Chen, Kemin Zhou and Jorge L. Aravena\\
Department of Electrical and Computer Engineering\\
Louisiana State University\\
Baton Rouge, LA 70803}

\date{Received in July 2002, Revised in February 2004}

 \maketitle

\begin{abstract}                          % Abstract of not more than 200 words.
The Clopper-Pearson confidence interval has ever been documented
as an exact approach in some statistics literature.  More
recently, such approach of interval estimation has been
introduced to probabilistic control theory and has been referred
as non-conservative in control community. In this note, we
clarify the fact that the so-called exact approach is actually
conservative. In particular, we derive analytic results
demonstrating the extent of conservatism in the context of
probabilistic robustness analysis. This investigation encourages
seeking better methods of confidence interval construction for
robust control purpose.
\end{abstract}

\section{Introduction}
Ever since  Stengel and Ray originated the concept of stochastic
robustness, there has been growing interest in developing
probabilistic methods for robust control. Significant
contributions have been made by a number of researchers (see,
e.g., \cite{BaTeFu:98,BarLag:97,BaLaTe:97,CaDaTe:00,ChZhAr:04},
\cite{KhaTik:96,PolSch:00,RaySte:93,RosBar:01,SteRay:91,TeBaDa:97,TemDab:99,VidBlo:01}
and the reference therein). A fundamental problem in the area of
probabilistic robustness analysis is to estimate the probability
that a certain robustness
 requirement is guaranteed for an uncertain dynamic system.
 The estimation of such probability relies essentially on Monte Carlo simulation.
 When an estimate of the probability is obtained from i.i.d. observations of fixed sample size,
 an important concern is how accurate this estimate is.
 To be useful, a numerical method must include a basis for error assessment.
 The Monte Carlo method is no exception.  Stengel and Ray \cite{RaySte:93,SteRay:91} first introduced
 the Clopper-Pearson confidence interval \cite{CloPea:34} to evaluate the accuracy of estimation
 in the context of robustness analysis.  In their works, such approach has been considered as non-conservative.

In this note, we would like to clarify the fact that the
Clopper-Pearson confidence interval is conservative.  The
erroneous understanding of the confidence interval is not due to
researchers in control area. Historically, such erroneous concept
can be traced back to some statistical literature. The
Clopper-Pearson confidence interval was usually referred as the
``exact'' confidence interval.  It has been documented as a
non-conservative approach in some statistics literature (see, for
example, Page 697-698 of \cite{Hald:52} and Page 95-103 of
\cite{Conover:77}). Although the so-called exact confidence
interval was proposed by Clopper and Pearson in 1934, its rigorous
probabilistic implication was not made clear until 1958 by
Cluniess-Ross \cite{Clunies-Ross:58}. Interestingly, it was
proved in \cite{Clunies-Ross:58} that such ``exact'' approach is
actually conservative.

Since the confidence coefficient is directly related to the risk
and safety of control systems in the context of probabilistic
design and analysis, a clear understanding of the conservatism of
confidence interval construction will help making the tradeoff
between the risk and performance enhancement. To this purpose, we
investigate the conservatism of the confidence interval.  We
obtain analytic results indicating that, in the scenario of rare
events (especially in the context of probabilistic robustness
analysis), the conservatism is not trivial and better methods of
confidence construction should be sought.

\section{Binomial Confidence Interval}

Let the probability space be denoted as $(\Omega, F,P)$ where
$\Omega, F, P$ are the sample space, the algebra of events and
the probability measure respectively. Let $X$ be a Bernoulli
random variable with distribution ${\rm Pr} \{ X=1\} =
\mathbb{P}_X, \;\;{\rm Pr} \{X=0\} = 1-\mathbb{P}_X$ where
$\mathbb{P}_X \in (0,1)$.   Let the sample size $N$ and confidence
parameter $\delta \in (0,1)$ be fixed. We refer an observation
with value 1 as a successful trial. Let $K$ denote the number of
successful trials during the $N$ i.i.d. sampling experiments. Let
$k = K(\omega)$ where $\omega$ is a sample point in the sample
space $\Omega$.

The classic Clopper-Pearson lower confidence limit
$L_{N,k,\delta}$ and
   upper confidence limit $U_{N,k,\delta}$ are given respectively by
\[
L_{N,k,\delta}:=\left\{\begin{array}{ll}
   0 \;\;\;&  {\rm if}\; k=0\\
   \underline{p} \; \;\;\;&
   {\rm if}\; k > 0
\end{array} \right.
\]
and \[
U_{N,k,\delta}:=\left\{\begin{array}{ll}
   1  &  {\rm if}\; k=N\\
   \overline{p} \;
   & {\rm if}\; k < N
\end{array} \right.
\]
where $\underline{p} \in (0,1)$ is the solution of equation
$\sum_{j=0}^{k-1} {N \choose j}
   \underline{p}^j (1- \underline{p})^{N-j} = 1-\frac{\delta}{2}$
   and $\overline{p} \in
(0,1)$ is the solution of equation $\sum_{j=0}^{k} {N \choose j}
   \overline{p}^j (1- \overline{p})^{N-j} = \frac{\delta}{2}$.
   Define random variable $L: \Omega \rightarrow [0,1]$ by
$L (\omega):=L_{N,K(\omega),\delta} \;\;\forall \omega \in \Omega$
and random variable $U: \Omega \rightarrow [0,1]$ by $U
(\omega):=U_{N,K(\omega),\delta} \;\;\forall \omega \in \Omega$.
Then the random interval $[L,U]$ is referred as the classic
Clopper-Pearson confidence interval. Its probabilistic
implication was quite often erroneously interpreted as
\[
{\rm Pr}\{ L \leq \mathbb{P}_X \leq U\} = 1-\delta
\]
or \[ {\rm Pr}\{ L < \mathbb{P}_X < U\} = 1-\delta.
\]
However, it was proved by Cluniess-Ross \cite{Clunies-Ross:58} in
1958 that
\[
{\rm Pr}\{ L \leq \mathbb{P}_X \leq U\} > 1-\delta \] and
\[ {\rm
Pr}\{ L < \mathbb{P}_X < U\} \geq 1-\delta.
\]
These inequalities have been demonstrated by numerical
experiments reported in the literature. For a better
understanding of the conservatism, especially in the context of
probabilistic robustness analysis, we shall investigate
analytically how conservative the Clopper-Pearson
interval can be.

\section{How Conservative?}

We refer the exact value of ${\rm Pr}\{ L \leq \mathbb{P}_X \leq
U\}$ as the {\it coverage probability}.  We have the following
results with regard to the conservatism of
 the Clopper-Pearson confidence interval.

\begin{theorem} \label{tight}

Let $N, \;\delta$ be fixed.  If $\mathbb{P}_X$ or $1-
\mathbb{P}_X$ is less than $1 - \left( \frac{\delta}{2}
\right)^{\frac{1}{N}}$, then the coverage probability will be at
least $1 - \frac{\delta}{2}$. Moreover, if $\left(
\frac{\delta}{2} \right)^{\frac{1}{N}} < \mathbb{P}_X < 1 - \left(
\frac{\delta}{2} \right)^{\frac{1}{N}}$, then the coverage
probability is $1$.
\end{theorem}

\bpf For the simplicity of notation, define $\mathcal{S} (N,k,x) :=
\sum_{j=0}^{k} {N \choose j} x^j (1- x)^{N-j}$ for $x \in (0,1)$.
Notice that $ \mathbb{P}_X \geq {\left( \frac{\delta}{2}
\right)}^{\frac{1}{N}}$ if and only if
\[
\mathcal{ S} (N,N-1,\mathbb{P}_X) = 1-{\mathbb{P}_X}^N
 \leq 1-\frac{\delta}{2}
 \]
and observe that,  for fixed $x \in (0,1)$, $\mathcal{ S}
(N,k,x)$ increases monotonically
 with respect to $k$, we have that
 $\mathcal{ S} (N,N-1,\mathbb{P}_X) \leq 1-\frac{\delta}{2}$ if and
only if
\[
\mathcal{ S} (N,k-1,\mathbb{P}_X) \leq
1-\frac{\delta}{2}\;\;\forall k \in\{1,\cdots,N\}. \] By Lemma
(3.8 a) on page 277 of \cite{Clunies-Ross:58}, we have that
$\mathcal{ S} (N,k-1,x)$ decreases monotonically with respect to
$x \in (0,1)$. Therefore, $ L_{N,k,\delta} < \mathbb{P}_X $ for
all $k \in\{1,\cdots,N\}$. Recall that $L_{N,N,\delta} =0 <
\mathbb{P}_X$, we thus have
\begin{equation} \label{eq1}
 \mathbb{P}_X \geq {\left( \frac{\delta}{2} \right)}^{\frac{1}{N}}
\;\Longleftrightarrow\;{\rm Pr} \{ L \leq \mathbb{P}_X \} =1.
\end{equation} Similarly, we can show that \begin{equation} \label{eq2} \mathbb{P}_X
\leq 1- {\left( \frac{\delta}{2} \right)}^{\frac{1}{N}}
\;\Longleftrightarrow\;{\rm Pr} \{ \mathbb{P}_X  \leq U\} =1.
\end{equation} For the case that $\mathbb{P}_X \geq {\left(
\frac{\delta}{2} \right)}^{\frac{1}{N}}$, by (\ref{eq1}) and
Bonferoni's inequality
\begin{eqnarray*}
{\rm Pr}\{ L \leq \mathbb{P}_X \leq U\} & \geq & {\rm Pr} \{ L
\leq
\mathbb{P}_X \} + {\rm Pr} \{ \mathbb{P}_X  \leq U\}-1\\
& = & {\rm Pr} \{ \mathbb{P}_X  \leq U\}\\
& > & 1- \frac{\delta}{2}.
\end{eqnarray*}
 Similarly, for the other case that $\mathbb{P}_X \leq 1- {\left( \frac{\delta}{2} \right)}^{\frac{1}{N}}$
 we can show that ${\rm Pr}\{ L \leq \mathbb{P}_X \leq U\} >
 1- \frac{\delta}{2}$.  Thus the first statement is proved.

 Finally, the proof of second statement is completed by making use of Bonferoni's
 inequality and inequalities~(\ref{eq1}) and ~(\ref{eq2}).

\epf

\begin{figure}
\begin{center}
\includegraphics[height=6.6cm]{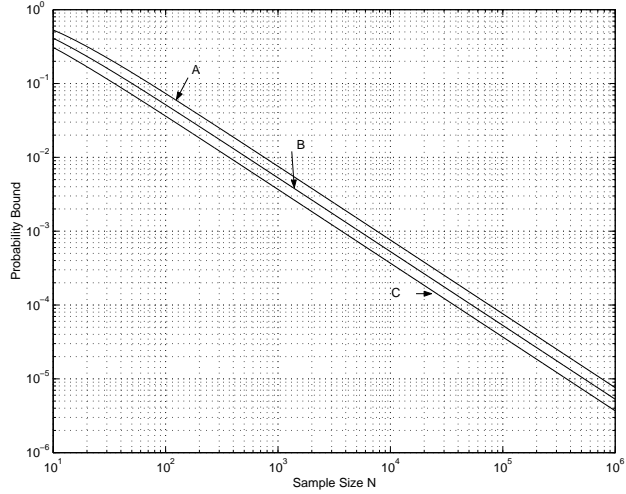}    % The printed column
\caption{ Conservatism of Clopper-Pearson Confidence Interval.
Probability Bound $ = 1 - \left( \frac{\delta}{2}
\right)^{\frac{1}{N}}$.  Plot A corresponds to $\delta = 0.001$,
plot B corresponds to $\delta = 0.01$, plot C
corresponds to $\delta = 0.05$.}  % width is 8.4 cm.
\label{fig_3}                                 % Size the figures
\end{center}                                 % accordingly.
\end{figure}

Figure~\ref{fig_3} shows bounds of binomial proportion for which
the true coverage probability of the Clopper-Pearson confidence
interval exceeds the prescribed confidence level by at least
$\frac{\delta}{2}$. It can be seen that the conservatism is
common in the scenario of studying rare events by the Monte Carlo
method. For example, when constructing a $99\%$ confidence
interval for the event of instability based on $10,000$ Monte
Carlo simulations, the true confidence level will be at least
$99.5\%$ if the probability of instability is smaller than
$0.0005$ (This number is the vertical coordinate of the point in
plot B with horizontal coordinate $10,000$).

\section{Conclusion}
It is demonstrated that the Clopper-Pearson confidence interval
is rather conservative for rare events. The coverage probability
can easily exceed the specified confidence level by at least
$\frac{\delta}{2}$ and can be $100\%$. Robustness issues such as
instability and performance violation are normally interpreted as
such rare events. Although the confidence parameter $\delta$ is
usually a small number, the impact can be enormous due to its
particular connection to the stability and performance of control
systems in the probabilistic robust control framework.  Our
investigation suggests seeking
better methods of confidence construction which are rigorous and less conservative.

\bibliographystyle{plain}        % Include this if you use bibtex
\bibliography{Chen}           % and a bib file to produce the
                                 % bibliography (preferred). The
                                 % correct style is generated by
                                 % Elsevier at the time of printing.

\end{document}